\documentclass[reqno]{amsart}
\usepackage{hyperref}
\usepackage{setspace}
\usepackage{graphicx}
\usepackage{amssymb}
\usepackage{mathrsfs}

\begin{document}
    \title[\hfilneg ]
    { Hilfer-Hadamard-type fractional differential equation with Cauchy-type problem }

 \author[\hfil\hfilneg]{Ahmad Y. A. Salamooni, D. D. Pawar }
 \address{Ahmad Y. A. Salamooni \newline
     School of Mathematical Sciences, Swami Ramanand Teerth Marathwada University, Nanded-431606, India}
      \email{ayousss83@gmail.com}

 \address{D. D. Pawar \newline
    School of Mathematical Sciences, Swami Ramanand Teerth Marathwada University, Nanded-431606, India}
     \email{dypawar@yahoo.com}


    \begin{abstract}
    In this paper, we consider the Cauchy-type problem (1.1) involving
    Hilfer-Hadamard-type fractional derivative for a nonlinear
    fractional differential equation. We prove an equivalence
    between the Cauchy-type problem (1.1) and Volterra integral
    equation(VIE), existence, and uniqueness. We present a slight
    generalization for the Gronwall inequality which was used in
    studying the continuous dependence of a solution for the Cauchy-type problem (1.2).
      \\\\ \textbf{AMS Classification- 34A08, 35R11}
    \end{abstract}

    \maketitle \numberwithin{equation}{section}
    \newtheorem{theorem}{Theorem}[section]
    \newtheorem{lemma}[theorem]{Lemma}
    \newtheorem{definition}[theorem]{Definition}
    \newtheorem{example}[theorem]{Example}

    \newtheorem{remark}[theorem]{Remark}
    \allowdisplaybreaks
  \[\textbf{1.Introduction}.\]
\par We consider the Cauchy-type problem
\begin{align*}
\quad\quad\quad\quad\quad\quad&~_{H}D^{\alpha,\beta}_{a+}x(t)=\varphi(t,x(t)),\quad\quad
n-1<\alpha<n,0\leq\beta\leq1,
\\&\quad\quad\quad\quad\quad\quad\quad\quad\quad\quad\quad\quad\quad\quad
\quad\quad\quad\quad\quad\quad\quad\quad\quad\quad\quad\quad\quad\quad\quad(1.1)\\&\quad
~_{H}D^{\gamma-j}_{a+}x(t)\big|_{t=a}=x_{a_{j}},\quad(j=1,2,...,n),\quad\quad\gamma=\alpha+\beta(n-\alpha),
\end{align*}
From the above initial condition and by definition 2.3(in this
paper), it is clear that
\[~_{H}D^{\gamma-j}_{a+}x(t)=\delta^{n-j}~_{H}I^{n-\gamma}_{a+}x(t).\] In
particular,
\begin{align*}
\quad\quad\quad\quad\quad\quad&~_{H}D^{\alpha,\beta}_{a+}x(t)=\varphi(t,x(t)),\quad\quad
0<\alpha<1,0\leq\beta\leq1,
\\&\quad\quad\quad\quad\quad\quad\quad\quad\quad\quad\quad\quad\quad\quad
\quad\quad\quad\quad\quad\quad\quad\quad\quad\quad\quad\quad\quad\quad\quad(1.2)\\&\quad
~_{H}I^{1-\gamma}_{a+}x(t)\big|_{t=a}=x_{a},\quad\quad\gamma=\alpha+\beta(1-\alpha),
\end{align*}
where $~_{H}D^{\alpha,\beta}_{a+}$ is the Hilfer-Hadamard-type
fractional derivative of order $\alpha$ and type $\beta,$[5,7].
Fractional differential equations have numerous applications in
science, physics, chemistry, and engineering,[1-4].
\par Recently, the theory and applications of fractional derivatives have received
considerable attention by researchers and authors. They have studied
the existence and uniqueness of solutions of fractional differential
equations on the different finite intervals such as the examples in
[5-8,13-15,17-21,23]and references therein. \par Some uses of the
Gronwall inequality with its applications to the fractional
derivatives and the continuous dependence for a solution on the
order of fractional differential equations under the initial
conditions are studied in [9,11,12,16,17].
\par In this paper, we found a variety of results for the initial
values problem (1.1), which are equivalent with (VIE), existence,
uniqueness, and Continuous dependence. In section 2, we present some
preliminaries. Section 3, contains the main results and is divided
into three parts. Part one dealt with the equivalence between the
Cauchy-type problem (1.1) and (VIE). Part two, which is section 3.1,
we proved the existence and uniqueness results for a solution of the
Cauchy-type problem (1.1) in the weighted space. The last part,
which is section 3.2, we found a slight generalization for the
Gronwall inequality and continuous dependence of the solution on the
order of the Cauchy-type (1.2) for Hadamard-type, and
Hilfer-Hadamard-type fractional differential equation under the
initial conditions.

\[\textbf{2.Preliminaries}\]\par In this section, we introduce some
notations, Lemmas, definitions and weighted spaces which are
important in developing some theories in this paper. For further
explanations, see [3]. \par Let $0<a<b<+\infty.$ Assume that
$C[a,b], AC[a,b],\quad and \quad C^{n}[a,b]$ be the spaces of
continuous, absolutely continuous, n-times continuous and
continuously differentiable functions on [a, b], respectively. And
let $L^{p}(a,b),with\quad p\geq1, $ be  the space of Lebesgue
integrable functions on (a, b). Moreover, we recall some of weighted
spaces[3] in definition2.1.
\\\textbf{\ Definition 2.1.[3]} Let $\Omega=[a,b]~(0<a<b<+\infty)$ is a finite interval and
$0\leq\mu<1,$ we introduce the weighted space $C_{\mu,\log}[a,b]$ of
continuous functions $\varphi$ on (a,b]
\[C_{\mu,\log}[a,b]=\{\varphi:(a,b]\rightarrow\mathbb{R}:[\log (t/a)]^{\mu}\varphi(t)\in C[a,b]\}\]
with the norm
\[\|\varphi\|_{C_{\mu,\log}}=\bigg\|[\log (t/a)]^{\mu}\varphi(t)\bigg\|_{C},\quad C_{0,\log}[a,b]=C[a,b].\]
And for $ n\in\mathbb{N}, and\quad \delta=t\frac{d}{dt}$
\[C^{n}_{\delta,\mu}[a,b]=\bigg\{\varphi:\|\varphi\|_{C^{n}_{\delta,\mu}}=\sum_{k=0}^{n-1}
\|\delta^{k}\varphi\|_{C}+\|\delta^{n}\varphi\|_{C_{\mu,\log}}\bigg\},C^{0}_{\delta,\mu}[a,b]=C_{\mu,\log}[a,b].\]
The space $C_{\mu,\log}[a,b]$ is the complete metric space defined
with the distance as
\[d(x_{1},x_{2})=\|x_{1}-x_{2}\|_{C_{\mu,\log}}[a,b]:=\max_{t\in[a,b]}\bigg|[\log (t/a)]^{\mu}\big[x_{1}(t)-x_{2}(t)\big]\bigg|\]
where $~\log(.)=\log_{e}(.)~$.
\\\textbf{\ Definition 2.2.[2,3]} Let $0<a<b<+\infty,$ the Hadamard fractional integral of order
$~\alpha\in \mathbb{R}^{+}~$for a function
$~\varphi:(a,\infty)\rightarrow\mathbb{R}~$ is defined
as\[_{H}I^{\alpha}_{a+}\varphi(t)=\frac{1}{\Gamma(\alpha)}_{a}\int^{t}(\log\frac{t}{\tau})^{\alpha-1}
\quad\frac{\varphi(\tau)}{\tau}d\tau,\quad\quad(t>a)\]
\\\textbf{\ Definition 2.3.[2,3]} Let $0<a<b<+\infty,$ the Hadamard fractional derivative of order
$~\alpha~$ applied to the function
$~\varphi:(a,\infty)\rightarrow\mathbb{R}$ is defined
as\[_{H}D^{\alpha}_{a+}\varphi(t)=\delta^{n}(_{H}I^{n-\alpha}_{a+}\varphi(t)),\quad
n-1<\alpha<n,\quad n=[\alpha]+1,\]
where$\quad~\delta^{n}=(t\frac{d}{dt})^{n}\quad~$and$~[\alpha]~$denotes
the integer part of the real number$~\alpha.~$
\\ \textbf{\ Lemma2.4.[3]} Let $n\in\mathbb{N}_{0}=\{0,1,2,...\}$ and let
$\mu_{1},\mu_{2}\in\mathbb{R},$ such that
$0\leq\mu_{1}\leq\mu_{2}<1.$
\\The following embeddings hold:
\[C_{\delta}^{n}[a,b] \longrightarrow C_{\delta,\mu_{1}}^{n}[a,b]\longrightarrow C_{\delta,\mu_{2}}^{n}[a,b], \]
with \[\|\varphi\|_{C_{\delta,\mu_{2}}^{n}}\leq
K_{\delta}\|\varphi\|_{C_{\delta,\mu_{1}}^{n}},\quad K_{\delta}=\min
\bigg[1,\bigg(\log (b/a)\bigg)^{\mu_{2}-\mu_{1}}\bigg],~~a\neq0\] In
particular,
\[C[a,b] \longrightarrow C_{\mu_{1},\log}[a,b]\longrightarrow C_{\mu_{2},\log}[a,b]\]
with \[\|\varphi\|_{C_{\mu_{2},\log}}\leq\bigg(\log
(b/a)\bigg)^{\mu_{2}-\mu_{1}}\|\varphi\|_{C_{\mu_{1},\log}},~~a\neq0\]
\\ \textbf{\ Lemma2.5.[3]}
\\$(a_{1})$If $\mathfrak{R}(\alpha)\geq0,\mathfrak{R}(\beta)\geq0,$ and
$~0<a<b<\infty,$ then
\begin{align*}
&\big[~_{H}I^{\alpha}_{a+}\big(\log(\tau/a)\big)^{\beta-1}\big](x)=\frac{\Gamma(\beta)}{\Gamma(\alpha+\beta)}(\log(t/a))^{\alpha+\beta-1},\quad
x>a\\&\big[~_{H}D^{\alpha}_{a+}\big(\log(\tau/a)\big)^{\beta-1}\big](x)=\frac{\Gamma(\beta)}{\Gamma(\alpha-\beta)}(\log(t/a))^{\alpha-\beta-1},\quad
x>a.
\end{align*}
\\$(a_{2})$ Let $\mathfrak{R}(\alpha)\geq0,n=[\mathfrak{R}(\alpha)]+1,$ and
$~0<a<b<\infty.$\par$~~$ The equality $(~_{H}D^{\alpha}_{a+}x)(t)=0$
is valid if, and only if,
\[x(t)=\sum_{k=1}^{n}c_{k}(\log(t/a))^{\alpha-k},\]
\par$~~$ where $c_{k}\in\mathbb{R}(k=1,2,...,n)$ are arbitrary constants.
\\$(a_{3})$ Let $\mathfrak{R}(\alpha)\geq0,\mathfrak{R}(\beta)\geq0,$
and $0\leq\mu<1.$ If $~0<a<b<\infty,$ then for $\varphi\in
C_{\mu,\log}[a,b]$
\[~_{H}I^{\alpha}_{a+}~_{H}I^{\beta}_{a+}\varphi=~_{H}I^{\alpha+\beta}_{a+}\varphi\]
\par$~~$ holds at any point $t\in (a,b].$ When $\varphi\in C[a,b],$ this
relation is valid at any point $t\in (a,b].$
\\\textbf{\ Theorem2.6.[3]} Let  $\mathfrak{R}(\alpha)\geq0,n=[\mathfrak{R}(\alpha)]+1,$ and
$~0<a<b<\infty.$ Also let $~_{H}I^{n-\alpha}_{a+}\varphi$ be the
Hadamard-type fractional integral of order $n-\alpha$ of the
function $\varphi.$ \par If $\varphi\in
C_{\mu,\log}[a,b]~(0\leq\mu<1)$ and
$~_{H}I^{n-\alpha}_{a+}\varphi\in  C_{\delta,\mu}^{n}[a,b],$ then
\[(_{H}I_{a+}^{\alpha}~_{H}D_{a+}^{\alpha}\varphi)(t)=
 \varphi(t)-\sum_{k=1}^{n}\frac{(\delta^{n-k}(_{H}I_{a+}^{n-\alpha}\varphi))(a)}{\Gamma(\alpha-k+1)}(\log\frac{t}{a})^{\alpha-k}\]
\\ \textbf{\ Lemma2.7.[3]} Let
$~0<a<b<\infty,\mathfrak{R}(\alpha)\geq0,n=[\mathfrak{R}(\alpha)]+1,$
and $0\leq\mathfrak{R}(\mu)<1.$ \par(a) If
$\mathfrak{R}(\mu)>\mathfrak{R}(\alpha)>0,$ then the fractional
integration operator $~_{H}I^{\alpha}_{a+}$ is  bounded from
$C_{\mu,\log}[a,b]$ into $C_{\mu-\alpha,\log}[a,b]$:
\[\|~_{H}I^{\alpha}_{a+}\varphi\|_{C_{\mu-\alpha,\log}}\leq k_{1}\|\varphi\|_{C_{\mu,\log}}\]
where \[k_{1}=\bigg(\log
(b/a)\bigg)^{\mathfrak{R}(\alpha)}\frac{\Gamma[\mathfrak{R}(\alpha)]|\Gamma(1-\mathfrak{R}(\mu))|}{|\Gamma(\alpha)|\Gamma(1+\mathfrak{R}(\alpha-\mu))}\]
In particular, $~_{H}I^{\alpha}_{a+}$ is  bounded in
$C_{\mu,\log}[a,b].$
\par(b) If $\mathfrak{R}(\mu)\leq\mathfrak{R}(\alpha),$ then the
fractional integration operator $~_{H}I^{\alpha}_{a+}$ is  bounded
from $C_{\mu,\log}[a,b]$ into $C[a,b]$:
\[\|~_{H}I^{\alpha}_{a+}\varphi\|_{C}\leq k_{2}\|\varphi\|_{C_{\mu,\log}}\]
where \[k_{2}=\bigg(\log
(b/a)\bigg)^{\mathfrak{R}(\alpha-\mu)}\frac{\Gamma[\mathfrak{R}(\alpha)]|\Gamma(1-\mathfrak{R}(\mu))|}{|\Gamma(\alpha)|\Gamma(1+\mathfrak{R}(\alpha-\mu))}\]
In particular, $~_{H}I^{\alpha}_{a+}$ is  bounded in
$C_{\mu,\log}[a,b].$
\\\textbf{\ Definition 2.8.[7]} Let$~~n-1<\alpha<n,~~0\leq\beta\leq 1,~~\varphi\in L^{1}(a,b).~$The
Hilfer-Hadamard fractional derivative $_{H}D^{\alpha,\beta}$ of
order$~\alpha~$ and type $~\beta~$of$~\varphi~$ is defined
as\[~(_{H}D^{\alpha,\beta}\varphi)(t)=\big(_{H}I^{\beta(n-\alpha)}(\delta)^{n}~_{H}I^{(n-\alpha)(1-\beta)}\varphi\big)(t)~\]
\[=\big(_{H}I^{\beta(n-\alpha)}(\delta)^{n}~_{H}I^{n-\gamma}\varphi\big)(t);\quad\gamma=\alpha+n\beta-\alpha\beta.\]
\[=\big(I^{\beta(n-\alpha)}_{H}D^{\gamma}\varphi\big)(t),\]
In particular, if $\quad0<\alpha<1,$ then
\[~(_{H}D^{\alpha,\beta}\varphi)(t)=\big(_{H}I^{\beta(1-\alpha)}\delta~_{H}I^{(1-\alpha)(1-\beta)}\varphi\big)(t)~\]
\[=\big(_{H}I^{\beta(1-\alpha)}\delta~_{H}I^{1-\gamma}\varphi\big)(t);\quad\gamma=\alpha+\beta-\alpha\beta.\]
\[=\big(_{H}I^{\beta(1-\alpha)}_{H}D^{\gamma}\varphi\big)(t).\]
Where $_{H}I^{(.)}$ and $~_{H}D^{(.)}~$is the Hadamard fractional
integral and derivative defined by (2.2) and (2.3), respectively.
\\\textbf{\ Definition 2.9.[3,18]} Assume that $\varphi(x,y)$ is defined on set $(a,b]\times
G,G\subset\mathbb{R}.$ A function $\varphi(x, y)$ satisfies
Lipschitz condition with respect to $y,$ if for all $x\in (a, b]$
and for all $y_{1}, y_{2} \in G,$
\[|\varphi(x,y_{1})-\varphi(x, y_{2})|\leq L|y_{1}- y_{2}|\] where $L>0$ is Lipschitz constant.
\\\textbf{\ Definition 2.10.[5,15]} Let $0<\alpha<1,0\leq\beta\leq 1,$ the weighted space
$C^{\alpha,\beta}_{1-\gamma}[a,b]$ is defined by
\[C^{\alpha,\beta}_{1-\gamma}[a,b]=\{\varphi\in C_{1-\gamma}[a,b]:D^{\alpha,\beta}_{a+}\varphi\in C_{1-\gamma}[a,b]\},\gamma=\alpha+\beta-\alpha\beta.\]
\\ \textbf{\ Lemma2.11.[10]} Let $0<a<b<\infty,\alpha>0,0\leq\mu<1,$ and $\varphi\in C_{\mu,\log}[a,b].$ If $\alpha>\mu,$ then $_{H}I_{a+}^{\alpha}\varphi$ is
continuous on $[a,b]$ and
\[_{H}I_{a+}^{\alpha}\varphi(a)=\lim_{t\rightarrow a^{+}}
~_{H}I_{a+}^{\alpha}\varphi(t)=0.\]
\\ \textbf{\ Lemma2.12.[7]} Let $\mathfrak{R}(\alpha)>0,0\leq\beta\leq1,\gamma=\alpha+n\beta-\alpha\beta,
n-1<\gamma\leq n,~n=[\mathfrak{R}(\alpha)]+1$ and $0<a<b<\infty.$ if
$\varphi\in L^{1}(a,b)$ and $(_{H}I_{a+}^{n-\gamma}\varphi)(t)\in
AC_{\delta}^{n}[a,b],$ then
\[_{H}I_{a+}^{\alpha}~(_{H}D_{a+}^{\alpha,\beta}\varphi)(t)=_{H}I_{a+}^{\gamma}~(_{H}D_{a+}^{\gamma}\varphi)(t)=
 \varphi(t)-\sum_{j=1}^{n}\frac{(\delta^{(n-j)}(_{H}I_{a+}^{n-\gamma}\varphi))(a)}
 {\Gamma(\gamma-j+1)}(\log\frac{t}{a})^{\gamma-j}\]
\\ \textbf{\ Lemma2.13.[18]} Let $0<a<b<\infty,0\leq\mu<1,\varphi\in C_{\mu,\log}[a,c]$ and $\varphi\in C_{\mu,\log}[c,b].$
Then, $\varphi\in C_{\mu,\log}[a,b]$ and
\[\|\varphi\|_{C_{\mu,\log}[a,b]}\leq\max\bigg\{\|\varphi\|_{C_{\mu,\log}[a,c]},\bigg(\log
(b/a)\bigg)^{\mu}\|\varphi\|_{C[c,b]}\bigg\}.\]
\\\textbf{\ Theorem2.14.[3]} Let $(\mathrm{U},d)$  be
a non-empty complete metric space, let $0\leq\omega<1,$ and let
$\mathrm{T}:\mathrm{U}\rightarrow\mathrm{U}$ be the map such that,
for every $u,v\in\mathrm{U},$ the relation
\[d(\mathrm{T}u,\mathrm{T}v)\leq\omega~d(u,v),\quad\quad0\leq\omega<1,\] holds. Then, the operator $\mathrm{T}$ has a
unique fixed point $u^{\ast}\in\mathrm{U}.$ \par Furthermore, if
$\mathrm{T}^{k}(k\in\mathbb{N})$ is the sequence of operators
defined by
\[\mathrm{T}^{1}=\mathrm{T},\quad\quad\mathrm{T}^{k}=\mathrm{T}\mathrm{T}^{k-1}\in\mathbb{N}\setminus\{1\},\] then,
for any $u_{0}\in\mathrm{U}$ the sequence
$\{\mathrm{T}^{k}u_{0}\}_{k=1}^{\infty} $ converges to the above
fixed point $u^{\ast}.$
\\\[\textbf{3.~Main Results }\]
\\\textbf{\ Definition 3.1.} Let $~~n-1<\alpha<n,~0\leq\beta\leq 1,\gamma=\alpha+\beta-\alpha\beta$ and $0\leq\mu<1,$ we consider the
underlying spaces defined by
\[C^{\alpha,\beta}_{\delta;n-\gamma,\mu}[a,b]=\{\varphi\in C_{n-\gamma,\log}[a,b]:
~_{H}D^{\alpha,\beta}_{a+}\varphi\in C_{\mu,\log}[a,b]\},\] and
\[C^{\gamma}_{n-\gamma,\log}[a,b]=\{\varphi\in C_{n-\gamma,\log}[a,b]:
~_{H}D^{\gamma}_{a+}\varphi\in C_{n-\gamma,\log}[a,b]\},\] where
$C_{n-\gamma,\log}[a,b]$ and $C_{\mu,\log}[a,b]$ are weighted spaces
of continuous functions on (a,b] defined by
\[C_{\gamma,\log}[a,b]=\big\{\varphi:(a,b]\rightarrow\mathbb{R}:\big(\log t/a\big)^{\gamma}\varphi(t)\in C[a,b]\big\}.\]
\par In the next theorem, we studied the equivalence between the Cauchy-type problem (1.1), and (VIE) of the second kind
\[\quad x(t)=\sum_{k=1}^{n}\frac{x_{a_{k}}}{\Gamma(\gamma-k+1)}(\log(t/a))^{\gamma-k}+\frac{1}{\Gamma(\alpha)}
_{a}\int^{t}(\log(t/\tau))^{\alpha-1}\varphi(\tau,x(\tau))\frac{d\tau}{\tau},~~
t>a\quad(3.1)\]
\\\textbf{\ Theorem 3.2.} Let $n-1<\alpha<n,0\leq\beta\leq1,\gamma=\alpha+\beta(n-\alpha),$
and assume that $\varphi(.,x(.))\in C_{\mu,\log}[a,b]$ where
$\varphi:(a,b]\times\mathbb{R}\rightarrow\mathbb{R}$ be a function
for any $x\in C_{\mu,\log}[a,b](n-\gamma\leq\mu<n-\beta(n-\alpha)).$
If $x\in C^{\gamma}_{n-\gamma,\log}[a,b],$ then $x$ satisfies (1.1)
if, and only if, $x$ satisfies the integral equation (3.1).
\\\\\textbf{\ Proof.} First part, we will Prove the necessity.
\\Assume that $x\in C^{\gamma}_{n-\gamma,\log}[a,b],$ is a solution
of (1.1). We prove that $x$ is a solution of (3.1) as follows: \\By
the definition 3.1 of $C^{\gamma}_{n-\gamma,\log}[a,b],$ Lemma
2.7(b), and definition 2.3, we have
\[~_{H}I^{n-\gamma}_{a+}x\in C[a,b],\quad\quad
~_{H}D^{\gamma}_{a+}x=\delta^{n}~_{H}I^{n-\gamma}_{a+}x\in
C_{n-\gamma,\log}[a,b].\] Thus by definition 2.1, we get
\[~_{H}I^{n-\gamma}_{a+}x\in C^{n}_{\delta,n-\gamma}[a,b].\]
Now, by applying Theorem 2.6, we obtain
\[\quad\quad\quad~_{H}I^{\gamma}_{a+}~_{H}D^{\gamma}_{a+}x(t)=x(t)-\sum_{k=1}^{n}\frac{(\delta^{n-k}(_{H}I_{a+}^{n-\gamma}\varphi))(a)}
{\Gamma(\gamma-k+1)}(\log\frac{t}{a})^{\gamma-k} \quad\quad
t\in(a,b], \quad\quad\quad\quad\quad\quad\] or
\[\quad\quad\quad~_{H}I^{\gamma}_{a+}~_{H}D^{\gamma}_{a+}x(t)=x(t)-\sum_{k=1}^{n}\frac{x_{a_{k}}}
{\Gamma(\gamma-k+1)}(\log\frac{t}{a})^{\gamma-k},~\quad t\in(a,b],
\quad\quad\quad\quad (3.2)\] where $x_{a_{k}}$ comes from the
initial condition of (1.1). By our hypothesis $\varphi(.,x(.))\in
C_{\mu,\log}[a,b],$ and since $x\in C_{n-\gamma,\log}[a,b]\subset
C_{\mu,\log}[a,b],$ and also by Lemma 2.7, we can see that the
integral $~_{H}I^{\alpha}_{a+}\varphi(.,x(.))\in
C_{\mu-\alpha,\log}[a,b]$ for $\mu>\alpha$ and
$~_{H}I^{\alpha}_{a+}\varphi(.,x(.))\in C[a,b]$ for $\mu\leq\alpha.$
It follows, by applying the operator $_{H}I^{\alpha}_{a+}$ to both
sides of the problem of Cauchy-type(1.1), and Lemma 2.12, that
\[\quad\quad~_{H}I^{\gamma}_{a+}~_{H}D^{\gamma}_{a+}x=~_{H}I^{\alpha}_{a+}~_{H}D^{\alpha,\beta}_{a+}x
=~_{H}I^{\alpha}_{a+}(~_{H}D^{\alpha,\beta}_{a+}x)=~_{H}I^{\alpha}_{a+}\varphi,\quad
in\quad(a,b]\quad\quad
(3.3)\] From (3.2) and (3.3), we get
\[\quad\quad x(t)=\sum_{k=1}^{n}\frac{x_{a_{k}}}
{\Gamma(\gamma-k+1)}(\log\frac{t}{a})^{\gamma-k}+
~_{H}I^{\alpha}_{a+}[\varphi(\tau,x(\tau))](t),\quad\quad
t\in(a,b]\quad\quad\quad
(3.4)\] which is the (VIE)(3.1)
\\\par Second part, we will Prove the sufficiency.
\\Assume that $x\in C^{\gamma}_{n-\gamma,\log}[a,b],$ satisfies (3.1),
that is written as(3.4), then$~_{H}D^{\gamma}_{a+}x$ exists
and$~_{H}D^{\gamma}_{a+}x\in C_{n-\gamma,\log}[a,b].$ Now by
applying the operator $~_{H}D^{\gamma}_{a+}$ to both sides of (3.4),
we get
\[~_{H}D^{\gamma}_{a+}x(t)=~_{H}D^{\gamma}_{a+}\bigg[\sum_{k=1}^{n}\frac{x_{a_{k}}}
{\Gamma(\gamma-k+1)}(\log\frac{t}{a})^{\gamma-k}+
~_{H}I^{\alpha}_{a+}[\varphi(\tau,x(\tau))](t)\bigg].\] By using
Lemma 2.5$(a_{2})~and~(a_{3}),$ and definition 2.3, we obtain
\begin{align*}\quad\quad\quad\quad\quad\quad\quad
~_{H}D^{\gamma}_{a+}x&=~_{H}D^{\gamma}_{a+}\big[_{H}I^{\alpha}_{a+}\varphi\big]\\&=
\delta^{n}(_{H}I^{n-\gamma}_{a+}~_{H}I^{\alpha}_{a+}\varphi)\\&=
\delta^{n}(_{H}I^{n-\beta(n-\alpha)}\varphi)\\&=
~_{H}D^{\beta(n-\alpha)}_{a+}\varphi\quad\quad\quad\quad\quad
\quad\quad\quad\quad\quad\quad\quad\quad\quad\quad\quad\quad\quad\quad(3.5)
\end{align*}
From (3.5), and the hypothesis $~_{H}D^{\gamma}_{a+}x\in
C_{n-\gamma,\log}[a,b],$ we have
\[\quad\quad\quad\quad\quad\quad\quad\quad\quad~_{H}D^{\beta(n-\alpha)}_{a+}\varphi\in
C_{n-\gamma,\log}[a,b]\quad\quad\quad\quad\quad\quad\quad\quad\quad\quad\quad\quad\quad\quad\quad(3.6)\]
Now, by applying $~_{H}I^{\beta(n-\alpha)}_{a+}$ to both sides of
(3.5), we obtain
\begin{align*}
&\quad\quad\quad\quad(~_{H}I^{\beta(n-\alpha)}_{a+}~_{H}D^{\gamma}_{a+}x)(t)=
(~_{H}I^{\beta(n-\alpha)}_{a+}~_{H}D^{\beta(n-\alpha)}_{a+}\varphi(\tau,x(\tau)))(t)\\&that~is,\\
&\quad\quad\quad\quad~_{H}I^{\beta(n-\alpha)}_{a+}\delta^{n}(_{H}I^{n-\gamma}_{a+}x)(t)=
(~_{H}I^{\beta(n-\alpha)}_{a+}~_{H}D^{\beta(n-\alpha)}_{a+}\varphi(\tau,x(\tau)))(t)
\quad\quad\quad\quad\quad(3.7)
\end{align*}
Since
\[\quad\quad\quad\delta^{n}(_{H}I^{n-\beta(n-\alpha)}_{a+}\varphi(t,x(t)))=~_{H}D^{\beta(n-\alpha)}_{a+}\varphi(.,x(.))\in
C_{n-\gamma,\log}[a,b], \quad\quad\quad\quad\quad(3.8)\]  and
$\gamma>\beta(n-\alpha)$ and by definition 2.1, we have
$_{H}I^{n-\beta(n-\alpha)}_{a+}\varphi\in
C_{\delta;n-\gamma}^{n}[a,b]$ (also that which is found in the first
part of this proof, or by Lemma 2.7(b) with $\mu<n-\beta(n-\alpha)$,
for a continuity of $_{H}I^{n-\beta(n-\alpha)}_{a+}\varphi$). Then,
Theorem 2.6, with definition 2.8 allows us to write
\[\quad\quad\quad_{H}D^{\alpha,\beta}_{a+}x(t)=\varphi(t,x(t))-\sum_{k=1}^{n}\frac{(\delta^{n-k}(_{H}I_{a+}^{n-\beta(n-\alpha)}\varphi))(a)}
{\Gamma(\beta(k-\alpha))}(\log\frac{t}{a})^{\beta(n-\alpha)-k},\quad\quad(3.9)\]
since, $\mu<n-\beta(n-\alpha),$ then it follows by Lemma 2.11, that
\[\big[_{H}I^{n-\beta(n-\alpha)}_{a+}\varphi)\big](a)=0\]
Therefore, we can write the relation (3.9) as
\[~_{H}D^{\alpha,\beta}_{a+}x(t)=\varphi(t,x(t)),\quad\quad t\in(a,b].\]
Finally, we will show that the initial condition of (1.1) also
holds. For that, we apply
$_{H}D^{\gamma-j}_{a+}=\delta^{n-j}~_{H}I^{n-\gamma}_{a+}(j=1,2,...,n)$
to both sides of (3.4), and by Lemma 2.5$(a_{1})~and~(a_{3}),$ we
obtain
\[\quad\quad\quad\quad\quad_{H}D^{\gamma-j}_{a+}x(t)=x_{a_{j}}+\big[\delta^{n-j}(_{H}I^{n-\beta(n-\alpha)}_{a+}\varphi(\tau,x(\tau)))\big](t)
\quad\quad\quad\quad\quad\quad\quad\quad\quad(3.10)\] Now, taking
the limit as $t\rightarrow a,$ in (3.10), we get
\[_{H}D^{\gamma-j}_{a+}x(t)\big|_{t=a}=x_{a_{j}},\quad\quad(j=1,2,...,n).\]
The proof of this theorem is complete.
\\\\\textbf{\ Remark 3.3.} For $0<\alpha<1,$ Theorem 3.2 is reduced to
Theorem 21(see[10]).
\[\textbf{3.1. Results of Existence and Uniqueness}\]
In this section, we will prove the existence and uniqueness results
for a solution of the Cauchy-type problem (1.1) in the weighted
space $C^{\alpha,\beta}_{n-\gamma,\log}[a,b],$ by using the Banach
fixed point theorem. For that, we need the following Lemma.
\\\\\textbf{\ Lemma 3.1.1.} If $\mu\in\mathbb{R}(0\leq\mu<1),$ then the Hadamard-type fractional
integral operator
$_{H}I^{\alpha}_{a+}~with~\alpha\in\mathbb{C}(\mathfrak{R}(\alpha)>0)$
is bounded from $C_{\mu,\log}[a,b] ~into~ C_{\mu,\log}[a,b] $ such
that,
\[\quad\quad\quad\quad\quad\quad\|~_{H}I^{\alpha}_{a+}\varphi\|_{C_{\mu,\log}[a,b]}\leq
\frac{\Gamma(1-\mu)}{\Gamma(1+\alpha-\mu)}
(\log(t/a))^{\alpha}\|\varphi\|_{C_{\mu,\log}[a,b]}.\quad\quad\quad\quad\quad(3.1.1)\]
\\\textbf{\ Proof.} By Lemma 2.7, the result of this Lemma follows.Now we will prove
the inequality(3.1.1). By definition 2.1 of the weighted space
$C_{\mu,\log}[a,b],$ we have
\begin{align*}
\|~_{H}I^{\alpha}_{a+}\varphi\|_{C_{\mu,\log}[a,b]}&=
\big\|(\log(t/a))^{\mu}~_{H}I^{\alpha}_{a+}\varphi\big\|_{C[a,b]}\\&\leq
\big\|\varphi\big\|_{C_{\mu,\log}[a,b]}\big\|_{H}I^{\alpha}_{a+}(\log(t/a))^{-\mu}\big\|_{C_{\mu,\log}[a,b]}
\end{align*}
Now, by using Lemma 2.5$(a_{1})(with~\beta~replaced~ by~ 1-\mu)$ we
obtain
\[\|~_{H}I^{\alpha}_{a+}\varphi\|_{C_{\mu,\log}[a,b]}\leq
\frac{\Gamma(1-\mu)}{\Gamma(1+\alpha-\mu)}
(\log(t/a))^{\alpha}\|\varphi\|_{C_{\mu,\log}[a,b]}.\] Hence, the
proof of this Lemma is complete.
\\\\\textbf{\ Theorem 3.1.2.} Let $n-1<\alpha<n,0\leq\beta\leq1,\gamma=\alpha+\beta(n-\alpha),$
and assume that $\varphi(.,x(.))\in C_{\mu,\log}[a,b]$ where
$\varphi:(a,b]\times\mathbb{R}\rightarrow\mathbb{R}$ be a function
for any $x\in C_{\mu,\log}[a,b](n-\gamma\leq\mu<n-\beta(n-\alpha)),$
and satisfies the Lipschitz condition given in definition 2.9 with
respect to $x.$ Then, there exists a unique solution $x(t)$ for the
Cauchy-type problem (1.1) in the weighted space
$C^{\alpha,\beta}_{\delta;n-\gamma,\mu}[a,b].$
\\\\\textbf{\ Proof.}
First, we will prove the existence of the unique solution $x(t)\in
C_{n-\gamma,\log}[a,b].$ According to Theorem 3.2, it is sufficient
to prove the existence of the unique solution $x(t)\in
C_{n-\gamma,\log}[a,b]$ to the nonlinear (VIE)(3.1), and that is
based on Theorem 2.14(Banach fixed point theorem). Since the
equation(3.1) makes sense in any interval $[a,t_{1}]\subset[a,b],$
then we choose $t_{1}\in(a,b]$ such that the following estimate
holds,
\[\quad\quad\quad\quad\quad\quad\quad
\omega_{1}:=L~\frac{\Gamma(\gamma-n+1)}{\Gamma(\alpha+\gamma-n+1)}(\log(t_{1}/a))^{\alpha}~<~1,
\quad\quad\quad\quad\quad\quad\quad\quad\quad\quad(3.1.2)\] where
$L>0$ is a Lipschitz constant. So we will prove the existence of the
unique solution $x(t)\in C_{n-\gamma,\log}[a,t_{1}]$ to the
equation(3.1) on the interval $(a,t_{1}].$ For this we  know that
the space $C_{n-\gamma,\log}[a,t_{1}]$ is a complete metric space
defined with the distance as
\[\quad\quad\quad~d(x_{1},x_{2})=\|x_{1}-x_{2}\|_{C_{n-\gamma,\log}}[a,t_{1}]:=\max_{t\in[a,t_{1}]}
\bigg|[\log
(t/a)]^{n-\gamma}\big[x_{1}(t)-x_{2}(t)\big]\bigg|.\quad\quad(3.1.3)\]
Rewrite equation(3.1) as the following:
\[\quad\quad\quad\quad\quad\quad\quad\quad\quad\quad\quad\quad\quad x(t)=(\mathrm{T}x)(t),
\quad\quad\quad\quad\quad\quad\quad\quad\quad\quad\quad\quad\quad\quad\quad\quad(3.1.4)\]
where $\mathrm{T}$ is the operator defined by
\[\quad\quad\quad\quad\quad\quad\quad\quad\quad (\mathrm{T}x)(t)=x_{0}(t)+\big[~_{H}I^{\alpha}_{a+}\varphi(\tau,x(\tau))\big](t),
\quad\quad\quad\quad\quad\quad\quad\quad\quad\quad(3.1.5)\] with
\[\quad\quad\quad\quad\quad\quad\quad\quad\quad\quad x_{0}(t)=\sum_{k=1}^{n}\frac{x_{a_{k}}}{\Gamma(\gamma-k+1)}(\log(t/a))^{\gamma-k},
\quad\quad\quad\quad\quad\quad\quad\quad\quad(3.1.6)\] Now, we claim
that $\mathrm{T}$ maps from
$C_{n-\gamma,\log}[a,t_{1}]~into~C_{n-\gamma,\log}[a,t_{1}].$ In
fact, it is clear from (3.1.6) that $x_{0}(t)\in
C_{n-\gamma,\log}[a,t_{1}].$ And since $\varphi(t,x(t))\in
C_{n-\gamma,\log}[a,t_{1}],$ then, by Lemma 2.7, and Lemma 3.1.1
$[with~ \mu=n-\gamma, b=t_{1}~ and~ \varphi(.)=\varphi(.,x(.))],$
the integral in the right-hand side of (3.1.5) relevant to
$C_{n-\gamma,\log}[a,t_{1}],$ and thus $(\mathrm{T}x)(t)\in
C_{n-\gamma,\log}[a,t_{1}].$
\par Next, we will prove that $\mathrm{T}$ is the contraction; that
is, we prove that the following estimate holds:
\[\quad\quad\quad\quad\quad\big\|\mathrm{T}x_{1}-\mathrm{T}x_{2}\big\|_{C_{n-\gamma,\log}[a,t_{1}]}\leq
\omega_{1}\big\|x_{1}-x_{2}\big\|_{C_{n-\gamma,\log}[a,t_{1}]},~0<\omega_{1}<1.\quad\quad\quad\quad(3.1.7)\]
By equations (3.1.5), (3.1.6), using the Lipschitz condition given
in definition 2.9, and applying the estimate (3.1.1)$[with~
\mu=n-\gamma, b=t_{1}~ and~
\varphi(t)=\varphi(t,x_{1}(t))-\varphi(t,x_{2}(t))],$ we get
\begin{align*}\quad
\big\|\mathrm{T}x_{1}-\mathrm{T}x_{2}\big\|_{C_{n-\gamma,\log}[a,t_{1}]}&\\&\leq
\big\|~_{H}I^{\alpha}_{a+}\big[|\varphi(t,x_{1}(t))-\varphi(t,x_{2}(t))|\big]\big\|_{C_{n-\gamma,\log}[a,t_{1}]}\\&\leq~L~
\big\|~_{H}I^{\alpha}_{a+}\big[|x_{1}(t))-x_{2}(t))|\big]\big\|_{C_{n-\gamma,\log}[a,t_{1}]}\\&\leq~
\omega_{1}\|x_{1}-x_{2}\|_{C_{n-\gamma,\log}[a,t_{1}]},\quad
\quad\quad\quad\quad \quad\quad\quad\quad\quad\quad\quad(3.1.8)
\end{align*}
which yields (3.1.7),$0<\omega_{1}<1.$ According to (3.1.2), and by
apply the Theorem 2.14(Banach fixed point theorem), we obtain a
unique solution $x^{\ast}\in C_{n-\gamma,\log}[a,t_{1}]$ to
(VIE)(3.1) on the interval $(a,t_{1}].$ \par This solution
$x^{\ast}$ is given from a limit of the convergent sequence
$(\mathrm{T}^{m}x_{0}^{\ast})(t):$
\[\quad\quad\quad\quad\quad\quad\quad\quad\quad\quad\lim_{m\rightarrow\infty}
\big\|\mathrm{T}^{m}x_{0}^{\ast}-x^{\ast}\big\|_{C_{n-\gamma,\log}[a,t_{1}]}=0,
\quad\quad\quad\quad\quad\quad\quad\quad\quad\quad\quad(3.1.9)\]
where $x_{0}^{\ast}$ is any function in $C_{n-\gamma,\log}[a,t_{1}]$
and
\begin{align*}\quad\quad\quad\quad\quad\quad
(\mathrm{T}^{m}x_{0}^{\ast})(t)&=(\mathrm{T}\mathrm{T}^{m-1}x_{0}^{\ast})(t)\\&=
x_{0}(t)+\big[~_{H}I^{\alpha}_{a+}\varphi(\tau,(\mathrm{T}^{m-1}x_{0}^{\ast})(\tau))\big](t),\quad
m\in\mathbb{N}.\quad\quad\quad\quad(3.1.10)
\end{align*}
Let us put $x_{0}^{\ast}(t)=x_{0}(t)~with~x_{0}(t)$ defined by
(3.1.6). \par If we indicate as
$x_{m}(t):=(\mathrm{T}^{m}x_{0}^{\ast})(t),$ then it is clear that
\[\quad\quad\quad\quad\quad\quad\quad\quad\quad\quad\lim_{m\rightarrow\infty}
\big\|x_{m}(t)-x^{\ast}\big\|_{C_{n-\gamma,\log}[a,t_{1}]}=0.
\quad\quad\quad\quad\quad\quad\quad\quad\quad\quad\quad\quad~(3.1.11)\]
Next, we consider the interval $[t_{1},b].$ From the (VIE)(3.1), we
have
\begin{align*}\quad\quad\quad\quad\quad
x(t)&=\sum_{k=1}^{n}\frac{x_{a_{k}}}{\Gamma(\gamma-k+1)}(\log(t/a))^{\gamma-k}\\&\quad\quad+\frac{1}{\Gamma(\alpha)}
_{a}\int^{t_{1}}(\log(t/\tau))^{\alpha-1}\varphi(\tau,x(\tau))\frac{d\tau}{\tau}\quad
\quad\quad\quad\quad\quad\quad\quad\quad\quad(3.1.12)\\&\quad\quad+\frac{1}{\Gamma(\alpha)}
_{t_{1}}\int^{t}(\log(t/\tau))^{\alpha-1}\varphi(\tau,x(\tau))\frac{d\tau}{\tau}\\&=x_{01}+\frac{1}{\Gamma(\alpha)}
_{t_{1}}\int^{t}(\log(t/\tau))^{\alpha-1}\varphi(\tau,x(\tau))\frac{d\tau}{\tau},
\end{align*}
where $x_{01}$ is defined by
\begin{align*}\quad\quad\quad\quad\quad\quad\quad
x_{01}&=\sum_{k=1}^{n}\frac{x_{a_{k}}}{\Gamma(\gamma-k+1)}(\log(t/a))^{\gamma-k}\quad\quad\quad\quad
\quad\quad\quad\quad\quad\quad\quad\quad\quad(3.1.13)\\&\quad\quad+\frac{1}{\Gamma(\alpha)}
_{a}\int^{t_{1}}(\log(t/\tau))^{\alpha-1}\varphi(\tau,x(\tau))\frac{d\tau}{\tau},
\end{align*}
and is the known function. We note that $x_{01}\in
C_{n-\gamma,\log}[t_{1},b].$ Now, we will prove the existence of the
unique solution $x(t)\in C_{n-\gamma,\log}[t_{1},b]$ to the
equation(3.1) on the interval $(t_{1},b].$ Also, here we use Theorem
2.14(Banach fixed point theorem) for the space
$C_{n-\gamma,\log}[t_{1},t_{2}]$ where
$t_{2}\in(t_{1},b](with~t_{2}=t_{1}+h_{1},~h_{1}>0,~t_{2}\leq b),$
satisfies
\[\quad\quad\quad\quad\quad\quad\quad
\omega_{2}:=L~\frac{\Gamma(\gamma-n+1)}{\Gamma(\alpha+\gamma-n+1)}(\log(t_{2}/t_{1}))^{\alpha}~<~1,
\quad\quad\quad\quad\quad\quad\quad\quad\quad\quad(3.1.14)\] the
space $C_{n-\gamma,\log}[t_{1},t_{2}]$ is a complete metric space
defined with the distance as
\[\quad\quad\quad~d(x_{1},x_{2})=\|x_{1}-x_{2}\|_{C_{n-\gamma,\log}}[t_{1},t_{2}]=\max_{t\in[t_{1},t_{2}]}
\bigg|[\log
(t/a)]^{n-\gamma}\big[x_{1}(t)-x_{2}(t)\big]\bigg|.\quad\quad(3.1.15)\]
Also we can rewrite equation(3.1.11) as the following:
\[\quad\quad\quad\quad\quad\quad\quad\quad\quad\quad\quad\quad\quad x(t)=(\mathrm{T}x)(t),
\quad\quad\quad\quad\quad\quad\quad\quad\quad\quad\quad\quad\quad\quad\quad\quad(3.1.16)\]
where $\mathrm{T}$ is the operator given by
\[\quad\quad\quad\quad\quad\quad\quad\quad\quad (\mathrm{T}x)(t)=x_{01}(t)+\big[~_{H}I^{\alpha}_{t_{1}+}\varphi(\tau,x(\tau))\big](t),
\quad\quad\quad\quad\quad\quad\quad\quad\quad\quad(3.1.17)\] As in
the beginning part of this proof, since $x_{01}(t)\in
C_{n-\gamma,\log}[t_{1},t_{2}],$ since $\varphi(t,x(t))\in
C_{n-\gamma,\log}[t_{1},t_{2}],$ then, by Lemma 2.7, and Lemma 3.1.1
$[with~ \mu=n-\gamma, b=t_{2}~ and~ \varphi(.)=\varphi(.,x(.))],$
the integral in the right-hand side of (3.1.16) also belongs to
$C_{n-\gamma,\log}[t_{1},t_{2}],$ and thus $(\mathrm{T}x)(t)\in
C_{n-\gamma,\log}[t_{1},t_{2}].$ \par Furthermore, using the
Lipschitz condition given in definition 2.9, and applying the
estimate(3.1.1)$[with~ \mu=n-\gamma, b=t_{2}~ and~
\varphi(t)=\varphi(t,x_{1}(t))-\varphi(t,x_{2}(t))],$ we get
\begin{align*}\quad\quad
\big\|\mathrm{T}x_{1}-\mathrm{T}x_{2}\big\|_{C_{n-\gamma,\log}[t_{1},t_{2}]}&=
\big\|~_{H}I^{\alpha}_{t_{1}+}\varphi(t,x_{1}(t))-~_{H}I^{\alpha}_{t_{1}+}\varphi(t,x_{2}(t))\big\|_{C_{n-\gamma,\log}[t_{1},t_{2}]}\\&\leq
\big\|~_{H}I^{\alpha}_{t_{1}+}\big[|\varphi(t,x_{1}(t))-\varphi(t,x_{2}(t))|\big]\big\|_{C_{n-\gamma,\log}[t_{1},t_{2}]}\\&\leq~L~
\big\|~_{H}I^{\alpha}_{t_{1}+}\big[|x_{1}(t))-x_{2}(t))|\big]\big\|_{C_{n-\gamma,\log}[t_{1},t_{2}]}\\&\leq~
\omega_{2}\|x_{1}-x_{2}\|_{C_{n-\gamma,\log}[t_{1},t_{2}]},\quad
\quad\quad\quad\quad \quad\quad\quad\quad\quad\quad(3.1.18)
\end{align*}
This, together with (3.1.13),$0<\omega_{2}<1,$ indicates that
$\mathrm{T}$ is a contraction and by applying the Theorem
2.14(Banach fixed point theorem), we obtain a unique solution
$x_{1}^{\ast}\in C_{n-\gamma,\log}[t_{1},t_{2}]$ to (VIE)(3.1) on
the interval $(t_{1},t_{2}].$ Moreover, this solution $x_{1}^{\ast}$
is given from a limit of the convergent sequence
$(\mathrm{T}^{m}x_{01}^{\ast})(t):$
\[\quad\quad\quad\quad\quad\quad\quad\quad\quad\quad\lim_{m\rightarrow\infty}
\big\|\mathrm{T}^{m}x_{01}^{\ast}-x_{1}^{\ast}\big\|_{C_{n-\gamma,\log}[t_{1},t_{2}]}=0,
\quad\quad\quad\quad\quad\quad\quad\quad\quad\quad\quad(3.1.19)\]
where $x_{01}^{\ast}$ is any function in
$C_{n-\gamma,\log}[t_{1},t_{2}],$ and again we can put
$x_{01}^{\ast}(t)=x_{01}(t)$ defined by (3.1.12). Hence,
\[\quad\quad\quad\quad\quad\quad\quad\quad\quad\quad\lim_{m\rightarrow\infty}
\big\|x_{m}(t)-x_{1}^{\ast}\big\|_{C_{n-\gamma,\log}[t_{1},t_{2}]}=0,
\quad\quad\quad\quad\quad\quad\quad\quad\quad\quad\quad\quad~(3.1.20)\]
where
\begin{align*}
\quad\quad\quad\quad\quad\quad~x_{m}(t)&=(\mathrm{T}^{m}x_{01}^{\ast})(t)\\&=
x_{01}(t)+\frac{1}{\Gamma(\alpha)}
_{t_{1}}\int^{t}(\log(t/\tau))^{\alpha-1}\varphi(\tau,x(\tau))\frac{d\tau}{\tau},\quad\quad\quad\quad
\quad\quad\quad(3.1.21)
\end{align*}
\par Next, if $t_{2}\neq b,$ we consider the interval
$[t_{2},t_{3}],$ such that
$t_{3}=t_{2}+h_{2},~with~h_{2}>0,~t_{3}\leq b$ and satisfies
\[\quad\quad\quad\quad\quad\quad\quad\quad
\omega_{3}:=L~\frac{\Gamma(\gamma-n+1)}{\Gamma(\alpha+\gamma-n+1)}(\log(t_{3}/t_{2}))^{\alpha}~<~1,
\quad\quad\quad\quad\quad\quad\quad\quad\quad(3.1.22)\] By using
same arguments as above, we conclude that there exists a unique
solution $x_{2}^{\ast}\in C_{n-\gamma,\log}[t_{2},t_{3}]$ to
(VIE)(3.1) on $[t_{2},t_{3}].$ If $t_{3}\neq b$, then we continue
the previous process until we get a unique solution $x(t)$ to the
(VIE)(3.1), and $x(t)=x_{i}^{\ast}$ such that $x_{i}^{\ast}\in
C_{n-\gamma,\log}[t_{i-1},t_{i}],~for~~i=1,2,...,L,$ where
$a=t_{0}<t_{1}<t_{2}<...<t_{L}=b,$ and
\[\quad\quad\quad\quad\quad\quad\quad\quad
\omega_{i+1}:=L~\frac{\Gamma(\gamma-n+1)}{\Gamma(\alpha+\gamma-n+1)}(\log(t_{i}/t_{i-1}))^{\alpha}~<~1,
\quad\quad\quad\quad\quad\quad\quad(3.1.23)\] Thus, by using Lemma
2.13, it yields that there exists a unique solution $x(t)\in
C_{n-\gamma,\log}[a,b]$ to the (VIE)(3.1), on the whole interval
[a,b].
\par Therefore, $x(t)\in
C_{n-\gamma,\log}[a,b]$ is a unique solution to the Cauchy-type
problem (1.1). \par Finally, we will show that such unique solution
$x(t)\in C_{n-\gamma,\log}[a,b]$ is in the weighted space
$C^{\alpha,\beta}_{n-\gamma,\mu}[a,b].$ By definition 3.1, it is
sufficient to prove that $~_{H}D^{\alpha,\beta}_{a+}x\in
C_{\mu,\log}[a,b].$ From the above proof, a solution $x(t)\in
C_{n-\gamma,\log}[a,b]$ is a limit of the sequence $x_{m}(t)\in
C_{n-\gamma,\log}[a,b]$ such that
\[\quad\quad\quad\quad\quad\quad\quad\quad\quad\quad\quad\lim_{m\rightarrow\infty}
\big\|x_{m}-x\big\|_{C_{n-\gamma,\log}[a,b]}=0,
\quad\quad\quad\quad\quad\quad\quad\quad\quad\quad\quad\quad~(3.1.24)\]
Hence, by using equation (1.1), Lipschitz condition given in
definition 2.9, and Lemma 2.4, we have
\begin{align*}
\quad\quad\quad\big\|~_{H}D^{\alpha,\beta}_{a+}x_{m}(t)-~_{H}D^{\alpha,\beta}_{a+}x(t)\big\|_{C_{\mu,\log}[a,b]}&=\big\|\varphi(t,x_{m}(t))-
\varphi(t,x(t))\big\|_{C_{\mu,\log}[a,b]}\\&\leq~L~(\log(b/a))^{\mu-n+\gamma}
\big\| x_{m}(t)-x(t)\big\|_{C_{n-\gamma,\log}[a,b]}\\&
\quad\quad\quad\quad\quad\quad\quad\quad\quad\quad\quad\quad\quad\quad\quad\quad(3.1.25)
\end{align*}
Clearly, the equations (3.1.24) and (3.1.25), yield that
\[\quad\quad\quad\quad\quad\quad\quad\quad\quad\lim_{m\rightarrow\infty}
\big\|~_{H}D^{\alpha,\beta}_{a+}x_{m}(t)-~_{H}D^{\alpha,\beta}_{a+}x(t)\big\|_{C_{\mu,\log}[a,b]}=0,
\quad\quad\quad\quad\quad\quad(3.1.26)\] and hence
$(~_{H}D^{\alpha,\beta}_{a+}x)\in~C_{\mu,\log}[a,b].$ Thus, the
proof of this theorem is complete.
\\\\\textbf{\ Remark 3.1.3.} For $0<\alpha<1,$ Theorem 3.1.2 reduced to
Theorem 22(see[10]).
\[\textbf{3.2. Results of Continuous Dependence }\]
In this section, firstly we wish to find a slight generalization for
the Gronwall inequality which can be used in the study of the
continuous dependence of a solution for the Cauchy-type problem
(1.2) of Hilfer-Hadamard-type fractional differential equation. And
the proof of the next Lemma is based on an iteration argument.
\\\\\textbf{\ Lemma 3.2.1.} Let $\alpha>~0,~u(t)$ be a nonnegative function locally integrable
on $a\leq t < T,(some~~a>~0,~~T\leq~+\infty)$ and $\psi(t)$ is a
nonnegative, nondecreasing continuous function defined on $a\leq t <
T,~\psi(t)\leq C$(constant), and Let $v(t)$ be nonnegative and
locally integrable on $a\leq t < T,$ with
\[\quad\quad\quad\quad\quad v(t)~\leq~u(t)+\psi(t)~_{a}\int^{t}(\log(t/\tau))^{\alpha-1}~\frac{v(\tau)}{\tau}d\tau,\quad\quad t\in[a,T)
\quad\quad\quad\quad\quad(3.2.1)\] Then
\[\quad v(t)~\leq~u(t) +~_{a}\int^{t}\Bigg[\sum_{k=1}^{\infty}\frac{\big(\psi(t)\Gamma(\alpha)\big)^{k}}{\Gamma(k\alpha)}(\log(t/\tau))^{k\alpha-1}
~\frac{v(\tau)}{\tau}\Bigg]d\tau,\quad\quad
t\in[a,T)\quad\quad(3.2.2)\]
\\\textbf{\ Proof.}
Assume that $M
\theta(t)=\psi(t)~_{a}\int^{t}(\log(t/\tau))^{\alpha-1}~\frac{\theta(\tau)}{\tau}d\tau,\quad
t\geq~a,$ for locally integrable functions $\theta.$ Then equation
(3.2.1), we can be written as
\[ v(t)~\leq~u(t)+M v(t)\] implies that
\[\quad\quad\quad\quad\quad\quad\quad\quad\quad\quad v(t)~\leq~\sum_{j=0}^{k-1}M^{j} u(t)+M^{k}v(t)
\quad\quad\quad\quad\quad\quad\quad\quad\quad\quad\quad\quad\quad\quad(3.2.3)\]
Now, we will prove that,
\[\quad\quad\quad\quad\quad\quad\quad\quad\quad\quad M^{k}v(t)~\leq~_{a}\int^{t}
\frac{\big(\psi(t)\Gamma(\alpha)\big)^{k}}{\Gamma(k\alpha)}(\log(t/\tau))^{k\alpha-1}
~\frac{v(\tau)}{\tau}d\tau,\quad\quad\quad\quad(3.2.4)\] and
$M^{k}v(t)\longrightarrow0,\quad as \quad
k\longrightarrow+\infty,\quad for ~each ~t\in[a,T),$ by using the
mathematical induction method. It is easy to show that relation
(3.2.4), true for $k=1.$ Assume that it is true for some $k=l.$ we
will show that it is true for $k=l+1,$ by hypothesis of induction
which yields that
\[M^{l+1}v(t)~\leq~(\psi(t))~~_{a}\int^{t}(\log(t/\tau))^{\alpha-1}\Bigg[~_{a}\int^{\tau}
\frac{\big(\psi(\tau)\Gamma(\alpha)\big)^{l}}{\Gamma(l\alpha)}(\log(\tau/s))^{l\alpha-1}
~\frac{v(s)}{s}ds\Bigg]d\tau.(3.2.5)\] Since $\psi(t)$ is a
nondecreasing, then the relation (3.2.5) is written as
\[\quad M^{l+1}v(t)~\leq~(\psi(t))^{l+1}~~_{a}\int^{t}(\log(t/\tau))^{\alpha-1}\Bigg[~_{a}\int^{\tau}
\frac{\big(\Gamma(\alpha)\big)^{l}}{\Gamma(l\alpha)}(\log(\tau/s))^{l\alpha-1}
~\frac{v(s)}{s}ds\Bigg]d\tau.(3.2.6)\] By interchanging the order of
integration, we can write the relation (3.2.6), as
\begin{align*}
M^{l+1}v(t)&\leq~(\psi(t))^{l+1}_{a}\int^{t}\Bigg[~_{s}\int^{t}
\frac{\big(\Gamma(\alpha)\big)^{l}}{\Gamma(l\alpha)}(\log(t/\tau))^{\alpha-1}
(\log(\tau/s))^{l\alpha-1}~\frac{d\tau}{\tau}\Bigg]\frac{v(s)}{s}ds\\&=~_{a}\int^{t}
\frac{\big(\psi(t)\Gamma(\alpha)\big)^{l+1}}{\Gamma((l+1)\alpha)}(\log(t/\tau))^{(l+1)\alpha-1}
~\frac{v(\tau)}{\tau}d\tau,\quad\quad\quad\quad\quad\quad\quad\quad\quad\quad\quad(3.2.7)
\end{align*}
where the integral,
\begin{align*}
_{s}\int^{t}(\log(t/\tau))^{\alpha-1}(\log(\tau/s))^{l\alpha-1}~\frac{d\tau}{\tau}&=
(\log(t/s))^{l\alpha+\alpha-1}~_{0}\int^{1}(1-z)^{\alpha-1}~
z^{l\alpha-1} dz \\&=(\log(t/s))^{(l+1)\alpha-1} B(l\alpha,\alpha)
\\&=(\log(t/s))^{(l+1)\alpha-1}\frac{\Gamma(\alpha)\Gamma(l\alpha)}{\Gamma((l+1)\alpha)},
\end{align*}
is done with the help of a substitution $\log\tau=\log s+z(\log
t-\log s),$ and a definition of beta function[1].Hence, the relation
(3.2.4) is proved.
\par Since $\psi(t)\leq C,$ and
\[M^{k}v(t)~\leq~_{a}\int^{t}
\frac{\big(\psi(t)\Gamma(\alpha)\big)^{k}}{\Gamma(k\alpha)}(\log(t/\tau))^{k\alpha-1}
~\frac{v(\tau)}{\tau}d\tau~\longrightarrow0~~as~~k\longrightarrow+\infty,~for~~a\leq
t<T,\]
Then the proof of this Lemma is complete.
\par Before studying the continuous dependence the
Cauchy-type problem (1.2), we will discuss some results for the
Cauchy-type problem of Hadamard fractional differential
equation[3,p.213], which given in the form
\begin{align*}
\quad\quad\quad\quad\quad\quad&~_{H}D^{\alpha}_{a+}x(t)=\varphi(t,x(t)),\quad\quad0<\alpha<1,
\\&\quad\quad\quad\quad\quad\quad\quad\quad\quad\quad\quad\quad\quad\quad\quad\quad
\quad\quad\quad\quad\quad\quad\quad\quad\quad\quad\quad\quad\quad\quad\quad(3.2.8)\\&\quad
~_{H}I^{1-\alpha}_{a+}x(t)\big|_{t=a}=x_{a},\quad\quad
\end{align*}
and for that we reduce (3.2.8) to the integral equation
\[\quad\quad\quad x(t)=\frac{x_{a}}{\Gamma(\alpha)}(\log(t/a))^{\alpha-1}+\frac{1}{\Gamma(\alpha)}
_{a}\int^{t}(\log(t/\tau))^{\alpha-1}\varphi(\tau,x(\tau))\frac{d\tau}{\tau},\quad
t>a\quad\quad\quad(3.2.9)\] where the equation (3.2.9) is the
equivalent to the initial value problem (3.2.8), see [3,p.213]. By
using the previous Lemma, we present a continuous dependence of the
solution on the order of the Cauchy-type problem of Hadamard-type
fractional differential equation.
\\\\\textbf{\ Theorem 3.2.2.} Let $0<\alpha-\delta\leq1,$ where $\alpha,\delta>0.$ Suppose that a
function $\varphi$ is continuous and satisfying the Lipschitz
condition (given in definition 2.9) in $\mathbb{R}.$ For $a\leq
t\leq h < b,$ Let $x$ be a solution of the initial value problem
(3.2.8), and $\tilde{x}$ be a solution of the initial value problem
\begin{align*}
\quad\quad\quad\quad\quad\quad&~_{H}D^{\alpha-\delta}_{a+}\tilde{x}(t)=\varphi(t,\tilde{x}(t)),\quad\quad0<\alpha<1,
\\&\quad\quad\quad\quad\quad\quad\quad\quad\quad\quad\quad\quad\quad\quad\quad\quad
\quad\quad\quad\quad\quad\quad\quad\quad\quad\quad\quad\quad\quad\quad\quad(3.2.10)\\&\quad
~_{H}I^{1-(\alpha-\delta)}_{a+}\tilde{x}(t)\big|_{t=a}=\tilde{x}_{a},\quad\quad
\end{align*}
Then, for $a< t\leq h,$ the estimate of the following
\[\big|\tilde{x}(t)-x(t)\big|\leq H(t)+~_{a}\int^{t}\Bigg[\sum_{k=1}^{\infty}\bigg(\frac{L\Gamma(\alpha-\delta)}{\Gamma(\alpha)}\bigg)^{k}
 ~\frac{(\log(t/\tau))^{k(\alpha-\delta)-1}}{\Gamma(k(\alpha-\delta))}\frac{H(\tau)}{\tau}\Bigg]d\tau,\]
hold, where
\begin{align*}
\quad\quad\quad\quad~H(t)&=\Bigg|\frac{\tilde{x}_{a}}{\Gamma((\alpha-\delta))}(\log(t/a))^{(\alpha-\delta)-1}-
\frac{x_{a}}{\Gamma(\alpha)}(\log(t/a))^{\alpha-1}\Bigg|\\&\quad+\|\varphi\|\Bigg|\frac{(\log(t/a))^{\alpha-\delta}}{\Gamma(\alpha-\delta+1)}-
\frac{(\log(t/a))^{\alpha-\delta}}{(\alpha-\delta)\Gamma(\alpha)}\Bigg|\\&\quad+\|\varphi\|\Bigg|
\frac{(\log(t/a))^{\alpha-\delta}}{(\alpha-\delta)\Gamma(\alpha)}-
\frac{(\log(t/a))^{\alpha}}{\Gamma(\alpha+1)}\Bigg|\quad\quad\quad\quad\quad\quad
\quad\quad\quad\quad\quad\quad\quad\quad\quad(3.2.11)
\end{align*}
and
\[\|\varphi\|=~\max_{a\leq~t\leq h}|\varphi(t,x(t))|.\]
\\\textbf{\ Proof.} The Solutions of the initial value problems (3.2.8), and (3.2.10)
are given by
\[\quad x(t)=\frac{x_{a}}{\Gamma(\alpha)}(\log(t/a))^{\alpha-1}+\frac{1}{\Gamma(\alpha)}
_{a}\int^{t}(\log(t/\tau))^{\alpha-1}\frac{\varphi(\tau,x(\tau))}{\tau}d\tau,\quad\]
and
\[\quad \tilde{x}(t)=\frac{x_{a}}{\Gamma(\alpha-\delta)}(\log(t/a))^{\alpha-\delta-1}+\frac{1}{\Gamma(\alpha-\delta)}
_{a}\int^{t}(\log(t/\tau))^{\alpha-\delta-1}\frac{\varphi(\tau,\tilde{x}(\tau))}{\tau}d\tau,\quad\]
It follows that
\begin{align*}
\quad\quad\big|\tilde{x}-x(t)\big|&=\Bigg|\frac{\tilde{x}_{a}}{\Gamma(\alpha-\delta)}(\log(t/a))^{\alpha-\delta-1}-
\frac{x_{a}}{\Gamma(\alpha)}(\log(t/a))^{\alpha-1}\\&\quad\quad+\frac{1}{\Gamma(\alpha-\delta)}
_{a}\int^{t}(\log(t/\tau))^{\alpha-\delta-1}\frac{\varphi(\tau,\tilde{x}(\tau))}{\tau}d\tau\\&\quad\quad\quad-\frac{1}{\Gamma(\alpha)}
_{a}\int^{t}(\log(t/\tau))^{\alpha-1}\frac{\varphi(\tau,x(\tau))}{\tau}d\tau\Bigg|\\&\leq
\Bigg|\frac{\tilde{x}_{a}}{\Gamma(\alpha-\delta)}(\log(t/a))^{\alpha-\delta-1}-
\frac{x_{a}}{\Gamma(\alpha)}(\log(t/a))^{\alpha-1}\Bigg|\\&\quad\quad+
\Bigg|_{a}\int^{t}\Bigg[\frac{(\log(t/\tau))^{\alpha-\delta-1}}{\Gamma(\alpha-\delta)}-
\frac{(\log(t/\tau))^{\alpha-\delta-1}}{\Gamma(\alpha)}\Bigg]\frac{|\varphi(\tau,\tilde{x}(\tau))|}{\tau}d\tau\Bigg|\\&\quad\quad+
\Bigg|\frac{1}{\Gamma(\alpha)}_{a}\int^{t}\bigg[(\log(t/\tau))^{\alpha-\delta-1}-(\log(t/\tau))^{\alpha-1}\bigg]
\frac{|\varphi(\tau,x(\tau))|}{\tau}d\tau\Bigg|\\&\quad\quad+\Bigg|_{a}\int^{t}\frac{(\log(t/\tau))^{\alpha-\delta-1}}{\Gamma(\alpha)}
\frac{|\varphi(\tau,\tilde{x}(\tau))-\varphi(\tau,x(\tau))|}{\tau}d\tau\Bigg|\\&\leq
\Bigg|\frac{\tilde{x}_{a}}{\Gamma(\alpha-\delta)}(\log(t/a))^{\alpha-\delta-1}-
\frac{x_{a}}{\Gamma(\alpha)}(\log(t/a))^{\alpha-1}\Bigg|\\&\quad\quad+\|\varphi\|
\Bigg|\frac{(\log(t/a))^{\alpha-\delta}}{\Gamma(\alpha-\delta+1)}-
\frac{(\log(t/a))^{\alpha-\delta}}{(\alpha-\delta)\Gamma(\alpha)}\Bigg|\\&\quad\quad+\|\varphi\|\Bigg|
\frac{(\log(t/a))^{\alpha-\delta}}{(\alpha-\delta)\Gamma(\alpha)}-
\frac{(\log(t/a))^{\alpha}}{\Gamma(\alpha+1)}\Bigg|\\&\quad\quad+\frac{L}{\Gamma(\alpha)}_{a}\int^{t}
(\log(t/\tau))^{\alpha-\delta-1}\frac{|\tilde{x}(\tau)-x(\tau)|}{\tau}d\tau
\end{align*}
Then, we have
\[\quad\quad\quad\quad\quad\quad\big|\tilde{x}-x(t)\big|\leq~H(t)+\frac{L}{\Gamma(\alpha)}_{a}\int^{t}
(\log(t/\tau))^{\alpha-\delta-1}\frac{|\tilde{x}(\tau)-x(\tau)|}{\tau}d\tau~\quad\quad\quad\quad\quad(3.2.12)\]
Where $H(t)$ is defined by (3.2.11). It follows by applying Lemma
3.2.1, that
\[\big|\tilde{x}-x(t)\big|\leq H(t)+~_{a}\int^{t}\Bigg[\sum_{k=1}^{\infty}\bigg(\frac{L\Gamma(\alpha-\delta)}{\Gamma(\alpha)}\bigg)^{k}
 ~\frac{(\log(t/\tau))^{k(\alpha-\delta)-1}}{\Gamma(k(\alpha-\delta))}\frac{H(\tau)}{\tau}\Bigg]d\tau,\]
Hence, the proof of theorem is complete.
\\ \par Next, we study the continuous dependence of the solution on the order of the
Cauchy-type problem of Hilfer-Hadamard-type fractional differential
equation by using the Lemma 3.2.1, for that we consider the initial
condition that given in (1.2), with $0<\alpha<1,~0\leq\beta\leq
1,~a\leq t < b,~(b\leq~+\infty),$ $\varphi :
[a,b)\times\mathbb{R}\rightarrow\mathbb{R}, $ and the solutions of
two initial value problems with a neighbouring orders and a
neighbouring initial values.
\\\\\textbf{\ Theorem 3.2.3.} Let $0<\alpha-\delta\leq1,$ where $\alpha,\delta>0.$ Suppose that a
function $\varphi$ is continuous and satisfying the Lipschitz
condition (given in definition 2.9) in $\mathbb{R}.$ For $a\leq
t\leq h < b,$ Let $x$ be a solution of the initial value problem
(1.2), and $\tilde{x}$ be a solution of the initial value problem
\begin{align*}
\quad\quad\quad\quad\quad\quad&~_{H}D^{\alpha-\delta,\beta}_{a+}\tilde{x}(t)=\varphi(t,\tilde{x}(t)),\quad\quad0<\alpha<1,~0\leq\beta\leq
1,\\&\quad\quad\quad\quad\quad\quad\quad\quad\quad\quad\quad\quad\quad\quad\quad\quad
\quad\quad\quad\quad\quad\quad\quad\quad\quad\quad\quad\quad\quad\quad\quad(3.2.13)\\&\quad
~_{H}I^{1-\gamma-\delta(\beta-1)}_{a+}\tilde{x}(t)\big|_{t=a}=\tilde{x}_{a},\quad\quad
\gamma=\alpha+\beta-\alpha\beta
\end{align*}
Then, for $a< t\leq h,$ the estimate of the following
\[\big|\tilde{x}(t)-x(t)\big|\leq F(t)+~_{a}\int^{t}\Bigg[\sum_{k=1}^{\infty}\bigg(\frac{L\Gamma(\alpha-\delta)}{\Gamma(\alpha)}\bigg)^{k}
 ~\frac{(\log(t/\tau))^{k(\alpha-\delta)-1}}{\Gamma(k(\alpha-\delta))}\frac{F(\tau)}{\tau}\Bigg]d\tau,\]
hold, where
\begin{align*}
\quad\quad\quad\quad~F(t)&=\Bigg|\frac{\tilde{x}_{a}}{\Gamma(\gamma+\delta(\beta-1))}(\log(t/a))^{\gamma+\delta(\beta-1)-1}-
\frac{x_{a}}{\Gamma(\gamma)}(\log(t/a))^{\gamma-1}\Bigg|\\&\quad+\|\varphi\|\Bigg|\frac{(\log(t/a))^{\alpha-\delta}}{\Gamma(\alpha-\delta+1)}-
\frac{(\log(t/a))^{\alpha-\delta}}{(\alpha-\delta)\Gamma(\alpha)}\Bigg|\\&\quad+\|\varphi\|\Bigg|
\frac{(\log(t/a))^{\alpha-\delta}}{(\alpha-\delta)\Gamma(\alpha)}-
\frac{(\log(t/a))^{\alpha}}{\Gamma(\alpha+1)}\Bigg|\quad\quad\quad\quad\quad\quad
\quad\quad\quad\quad\quad\quad\quad\quad\quad(3.2.14)
\end{align*}
and
\[\|\varphi\|=~\max_{a\leq~t\leq h}|\varphi(t,x(t))|.\]
\\\textbf{\ Proof.} By applying theorem 3.2, for the solutions $x,$ and $\tilde{x}$ of
the initial value problems (1.2), and (3.2.13) respectively. We can
write the solutions as
\[\quad x(t)=\frac{x_{a}}{\Gamma(\gamma)}(\log(t/a))^{\gamma-1}+\frac{1}{\Gamma(\alpha)}
_{a}\int^{t}(\log(t/\tau))^{\alpha-1}\frac{\varphi(\tau,x(\tau))}{\tau}d\tau,\quad\]
and
\[\quad \tilde{x}(t)=\frac{x_{a}}{\Gamma(\gamma+\delta(\beta-1))}(\log(t/a))^{\gamma+\delta(\beta-1)-1}+\frac{1}{\Gamma(\alpha-\delta)}
_{a}\int^{t}(\log(t/\tau))^{\alpha-\delta-1}\frac{\varphi(\tau,\tilde{x}(\tau))}{\tau}d\tau,\quad\]
It follows that
\begin{align*}
\quad\quad\big|\tilde{x}-x(t)\big|&=\Bigg|\frac{\tilde{x}_{a}}{\Gamma(\gamma+\delta(\beta-1))}(\log(t/a))^{\gamma+\delta(\beta-1)-1}-
\frac{x_{a}}{\Gamma(\gamma)}(\log(t/a))^{\gamma-1}\\&\quad\quad+\frac{1}{\Gamma(\alpha-\delta)}
_{a}\int^{t}(\log(t/\tau))^{\alpha-\delta-1}\frac{\varphi(\tau,\tilde{x}(\tau))}{\tau}d\tau\\&\quad\quad\quad-\frac{1}{\Gamma(\alpha)}
_{a}\int^{t}(\log(t/\tau))^{\alpha-1}\frac{\varphi(\tau,x(\tau))}{\tau}d\tau\Bigg|\\&\leq
\Bigg|\frac{\tilde{x}_{a}}{\Gamma(\gamma+\delta(\beta-1))}(\log(t/a))^{\gamma+\delta(\beta-1)-1}-
\frac{x_{a}}{\Gamma(\gamma)}(\log(t/a))^{\gamma-1}\Bigg|\\&\quad\quad+
\Bigg|_{a}\int^{t}\Bigg[\frac{(\log(t/\tau))^{\alpha-\delta-1}}{\Gamma(\alpha-\delta)}-
\frac{(\log(t/\tau))^{\alpha-\delta-1}}{\Gamma(\alpha)}\Bigg]\frac{|\varphi(\tau,\tilde{x}(\tau))|}{\tau}d\tau\Bigg|\\&\quad\quad+
\Bigg|\frac{1}{\Gamma(\alpha)}_{a}\int^{t}\bigg[(\log(t/\tau))^{\alpha-\delta-1}-(\log(t/\tau))^{\alpha-1}\bigg]
\frac{|\varphi(\tau,x(\tau))|}{\tau}d\tau\Bigg|\\&\quad\quad+\Bigg|_{a}\int^{t}\frac{(\log(t/\tau))^{\alpha-\delta-1}}{\Gamma(\alpha)}
\frac{|\varphi(\tau,\tilde{x}(\tau))-\varphi(\tau,x(\tau))|}{\tau}d\tau\Bigg|\\&\leq
\Bigg|\frac{\tilde{x}_{a}}{\Gamma(\gamma+\delta(\beta-1))}(\log(t/a))^{\gamma+\delta(\beta-1)-1}-
\frac{x_{a}}{\Gamma(\gamma)}(\log(t/a))^{\gamma-1}\Bigg|\\&\quad\quad+\|\varphi\|
\Bigg|\frac{(\log(t/a))^{\alpha-\delta}}{\Gamma(\alpha-\delta+1)}-
\frac{(\log(t/a))^{\alpha-\delta}}{(\alpha-\delta)\Gamma(\alpha)}\Bigg|\\&\quad\quad+\|\varphi\|\Bigg|
\frac{(\log(t/a))^{\alpha-\delta}}{(\alpha-\delta)\Gamma(\alpha)}-
\frac{(\log(t/a))^{\alpha}}{\Gamma(\alpha+1)}\Bigg|\\&\quad\quad+\frac{L}{\Gamma(\alpha)}_{a}\int^{t}
(\log(t/\tau))^{\alpha-\delta-1}\frac{|\tilde{x}(\tau)-x(\tau)|}{\tau}d\tau
\end{align*}
Then, we have
\[\quad\quad\quad\quad\quad\quad\big|\tilde{x}-x(t)\big|\leq~F(t)+\frac{L}{\Gamma(\alpha)}_{a}\int^{t}
(\log(t/\tau))^{\alpha-\delta-1}\frac{|\tilde{x}(\tau)-x(\tau)|}{\tau}d\tau~\quad\quad\quad\quad\quad(3.2.15)\]
Where $F(t)$ is defined by (3.2.14). It follows by applying Lemma
3.2.1, that
\[\big|\tilde{x}-x(t)\big|\leq F(t)+~_{a}\int^{t}\Bigg[\sum_{k=1}^{\infty}\bigg(\frac{L\Gamma(\alpha-\delta)}{\Gamma(\alpha)}\bigg)^{k}
 ~\frac{(\log(t/\tau))^{k(\alpha-\delta)-1}}{\Gamma(k(\alpha-\delta))}\frac{F(\tau)}{\tau}\Bigg]d\tau,\]
Hence, the proof of the theorem is complete.
\par In the next theorem, we shall make a small change of the
initial value condition that given in (1.2), as the following
\begin{align*}
\quad\quad\quad\quad\quad\quad&~_{H}D^{\alpha-\delta,\beta}_{a+}x(t)=\varphi(t,x(t)),\quad\quad0<\alpha<1,~0\leq\beta\leq
1,\\&\quad\quad\quad\quad\quad\quad\quad\quad\quad\quad\quad\quad\quad\quad\quad\quad
\quad\quad\quad\quad\quad\quad\quad\quad\quad\quad\quad\quad\quad\quad\quad(3.2.16)\\&\quad
~_{H}I^{1-\gamma}_{a+}x(t)\big|_{t=a}=x_{a}+\epsilon,\quad\quad
\gamma=\alpha+\beta-\alpha\beta,
\end{align*}
where $\epsilon$ be an arbitrary constant.
\\\\\textbf{\ Theorem 3.2.4.} Assume that suppositions of Theorem 3.2.2, hold. And assume $x(t)$ and
$\tilde{x}(t)$ are the solutions of the initial value problems (1.2)
and (3.2.16) respectively. Then the estimate of the following
\[\quad\quad\quad\quad\quad\big|\tilde{x}-x(t)\big|\leq~|\epsilon|~(\log(t/a))^{\gamma-1}
E_{\alpha,\gamma}\big[L(\log(t/a))^{\alpha}\big],\quad\quad
t\in(a,b]\quad\quad\quad\quad\quad(3.2.17)\] holds, where
$E_{\alpha,\gamma}(y)=\sum_{i=0}^{\infty}\frac{y^{i}}{\Gamma(n\alpha+\gamma)}$
be the  function of Mittag-Leffler, see[22].
\\\\\textbf{\ Proof.} According to the Theorem 3.1.2, we have
$x(t)=\lim_{n\rightarrow\infty}x_{n}(t),$ with $x_{0}(t),$ and
$x_{n}(t)$ are as defined in (3.1.3) and (3.1.4) respectively.
Clearly, for (3.2.16), we can write
$\tilde{x}=\lim_{n\rightarrow\infty}\tilde{x}_{n}(t),$ and
\[\quad\quad\quad\quad\quad\quad\quad\quad\quad\quad\quad \tilde{x}_{0}(t)=\frac{(x_{a}+\epsilon)}{\Gamma(\gamma)}~(\log(t/a))^{\gamma-1},\quad\quad\quad
\quad\quad\quad\quad\quad\quad\quad\quad\quad\quad(3.2.18)\]
\[\quad\quad\quad \tilde{x}_{n}(t)=\tilde{x}_{0}(t)+\frac{1}{\Gamma(\alpha)}
_{a}\int^{t}(\log(t/\tau))^{\alpha-1}\varphi(\tau,\tilde{x}_{n-1}(\tau))\frac{d\tau}{\tau},\quad
n\in\mathbb{N}\quad\quad\quad\quad\quad\quad\quad(3.2.19)\] It
follows from (3.1.3) and (3.2.18) that
\begin{align*}
\quad\quad\quad\quad\quad\quad\quad\quad&\big|x_{0}(t)-\tilde{x_{0}}(t)\big|=
\bigg|\frac{x_{a}}{\Gamma(\gamma)}(\log(t/a))^{\gamma-1}-\frac{(x_{a}+\epsilon)}{\Gamma(\gamma)}(\log(t/a))^{\gamma-1}\bigg|
\\&\big|x_{0}(t)-\tilde{x_{0}}(t)\big|\leq~|\epsilon|\frac{(\log(t/a))^{\gamma-1}}{\Gamma(\gamma)}
\quad\quad\quad\quad\quad\quad\quad\quad\quad\quad\quad\quad\quad\quad\quad(3.2.20)
\end{align*}
Now, by using equations (3.1.4) and (3.2.19), and applying the
Lipschitz condition given in definition 2.9, and the estimate
(3.2.20), we get
\begin{align*}
\quad\quad\quad\big|x_{1}(t)-\tilde{x_{1}}(t)\big|&=
\bigg|\epsilon\frac{(\log(t/a))^{\gamma-1}}{\Gamma(\gamma)}+\frac{1}{\Gamma(\alpha)}
_{a}\int^{t}(\log(t/\tau))^{\alpha-1}\frac{\big[\varphi(\tau,x_{0}(\tau))-\varphi(\tau,\tilde{x_{0}}(\tau))\big]}{\tau}d\tau\bigg|\\&\leq
|\epsilon|\frac{(\log(t/a))^{\gamma-1}}{\Gamma(\gamma)}+\frac{L}{\Gamma(\alpha)}_{a}\int^{t}
(\log(t/\tau))^{\alpha-1}\frac{|x_{0}(\tau)-\tilde{x_{0}}(\tau)|}{\tau}d\tau\\&\leq
|\epsilon|(\log(t/a))^{\gamma-1}\bigg[\frac{1}{\Gamma(\gamma)}+\frac{L(\log(t/a))^{\alpha}}{\Gamma(\alpha+\gamma)}\bigg]
\end{align*}
Then, we have
\[\quad\quad\quad\quad\quad\quad\quad\quad\big|x_{1}(t)-\tilde{x_{1}}(t)\big|\leq~|\epsilon|(\log(t/a))^{\gamma-1}
\sum_{i=0}^{1}\frac{L^{i}(\log(t/a))^{\alpha~i}}{\Gamma(\alpha~i+\gamma)}
\quad\quad\quad\quad\quad\quad\quad(3.2.21)\] Similarly, by using
the relation (3.2.21), it directly yields
\begin{align*}
\quad\quad\quad\big|x_{2}(t)-\tilde{x_{2}}(t)\big|&\leq
|\epsilon|\frac{(\log(t/a))^{\gamma-1}}{\Gamma(\gamma)}+\frac{L}{\Gamma(\alpha)}_{a}\int^{t}
(\log(t/\tau))^{\alpha-1}\frac{|x_{1}(\tau)-\tilde{x_{1}}(\tau)|}{\tau}d\tau\\&\leq
|\epsilon|(\log(t/a))^{\gamma-1}
\sum_{i=0}^{2}\frac{L^{i}(\log(t/a))^{\alpha~i}}{\Gamma(\alpha~i+\gamma)}
\quad\quad\quad\quad\quad\quad\quad\quad\quad\quad\quad\quad(3.2.22)
\end{align*}
And by using the mathematical induction method, we conclude that
\[\quad\quad\quad\quad\quad\quad\quad\quad\big|x_{n}(t)-\tilde{x_{n}}(t)\big|\leq~|\epsilon|(\log(t/a))^{\gamma-1}
\sum_{i=0}^{n}\frac{L^{i}(\log(t/a))^{\alpha~i}}{\Gamma(\alpha~i+\gamma)}
\quad\quad\quad\quad\quad\quad\quad(3.2.23)\] is taking a limit of
the summation as $n\longrightarrow\infty$ in relation (3.2.23).
Hence, we obtain
\begin{align*}
\quad\quad\quad\big|x_{n}(t)-\tilde{x_{n}}(t)\big|&\leq
|\epsilon|(\log(t/a))^{\gamma-1}
\sum_{i=0}^{\infty}\frac{L^{i}(\log(t/a))^{\alpha~i}}{\Gamma(\alpha~i+\gamma)}\\&=
E_{\alpha,\gamma}\big[L(\log(t/a))^{\alpha}\big].
\end{align*}
Then, the proof of this theorem is complete.


\begin{thebibliography}{99}


\bibitem{Po99} I. Podlubny;
\newblock{\it Fractional Differential Equations.}
\newblock{Academic Press, San Diego}, 1999

\bibitem{Sa-Ki-Ma93} S. G. Samko, A. A. Kilbas; O. I. Marichev;
\newblock{\it Fractional Integrals and Derivatives: Theory and Applications}.
\newblock{Gordon and Breach, New York (1993)}.
[Translation from the Russian edition, Nauka i Tekhnika, Minsk
(1987)]

\bibitem{Ki-Sr-Tr06} A. A. Kilbas, H. M. Srivastava, J. J. Trujillo;
\newblock{\it Theory and Applications of Fractional Differential Equations}.
\newblock{ North-Holland Mathematics Studies
204, Editor: Jan Van Mill, Elsevier, Amsterdam, The Netherlands, }
(2006).

\bibitem{Hi00}  Hilfer R.,
\newblock{\it Applications of Fractional Calculus in Physics}.
\newblock{World Scientific Publ. Co., Singapore}, 2000.

\bibitem{Qa-Fu-Sb12}  Qassim M.D., Furati K.M., Tatar N-e.
\newblock{\it On a differential equation involving Hilfer-
-hadamard fractional derivative.}.
  Abstract Appl Anal 2012;2012:17. Article ID 391062

\bibitem{Hi-Rp08}  Hilfer R.,
\newblock{\it . Threefold introduction to fractional derivatives.
In: Anomalous transport: foundations and applications;}.
  2008. p. 17-73

\bibitem{AS-DD} Ahmad Y. A. Salamooni, D. D. Pawar,
\newblock{\it Existence and uniqueness of boundary value
problems for Hilfer-Hadamard-type fractional differential
equations.}
 arXiv:1801.10400v1[math.AP] 31 Jan 2018.

\bibitem{JA-DD12} J. A. Nanware, D. B. Dhaigude,
\newblock{\it Existence and uniqueness of solution
of Riemann-Liouville fractional differential equations with integral
boundary conditions,}
\newblock{   Int. J. Nonlinear Sci., 14, (2012), 410-415.}

\bibitem{O.L}  O. Lipovan,
\newblock{\it A retarded Gronwall-like inequality and its applications.}
\newblock{J. Math. Anal. Appl. 252 (2000) 389-401.}

\bibitem{Ka013} M. D. Kassim and N. E. Tatar,
\newblock{\it Well-Posedness and Stability for a Differential Problem with
Hilfer-Hadamard Fractional Derivative,}
\newblock{Abst. Appl. Anal., Vol. 2013, (2013), 1–12. Article ID 605029.}

\bibitem{R.A-S.D-W.Z}  R.P. Agarwal, S. Deng, W. Zhang,
\newblock{\it Generalization of a retarded Gronwall-like inequality and its applications, }
\newblock{Appl. Math. Comput. 165 (2005) 599-612.}

\bibitem{BGP98}   B.G. Pachpatte,
\newblock{\it Inequalities for Differential and Integral Equations,  }
\newblock{in: Mathematics in Science and Engineering, vol. 197, Acad. Press, 1998.}

\bibitem{HR-Y-Z009}   Hilfer R., Y. Luchko, $\check{Z}.$ Tomovski,
\newblock{\it Operational method for solution of the
fractional differential equations with the generalized
Riemann-Liouville fractional derivatives,}
\newblock{Fract. Cal. Appl. Anal., (12), (2009), 299-318.}

\bibitem{KD-NJ002}  K. Diethelm, N. J. Ford,
\newblock{\it Analysis of fractional differential equations,}
\newblock{J. Math. Anal. Appl., Vol. 265, (2002), 229-248.}

\bibitem{KM-MK-N012}  K. M. Furati, M. D. Kassim and N.e.Tatar,
\newblock{\it Existence and uniqueness for
a problem involving Hilfer fractional derivative,}
\newblock{Computers Math. Appl.,(2012), 1616-1626.}

\bibitem{HY-JG-D012}  H. Ye, J. Gao, Y. Ding,
\newblock{\it A generalized Gronwall inequlity and its application to a fractional differential equation,}
\newblock{J. Math. Anal. Appl., 328, (2007),
1075-1081.}

\bibitem{DB-SB017} D. B. Dhaigude, *Sandeep P. Bhairat,
\newblock{\it Existence and uniqueness of solution of Cauchy-type problem for
Hilfer fractional differential equations,}
\newblock{arXiv:1704.02174v1 [Math.CA] 7 Apr 2017.
1075–1081.}

\bibitem{C-JL-Y010} C. Kou, J. Liu, and Y. Ye,
\newblock{\it Existence and uniqueness of solutions for the
Cachy-type problems of fractional differential equaitions, }
\newblock{Discrete Dyn. Nat. Soc., Article ID 142175, (2010), 1-15.}

\bibitem{HMS-ZT009} H.M. Srivastava, $\check{Z}.$ Tomovski,
\newblock{\it Fractional calculus with an integral operator containing a generalized Mittag-Leffler function in the kernel, }
\newblock{Applied Mathematics and Computation 211 (2009) 198-210.}

\bibitem{ZT-RH-HS10}  $\check{Z}$ivorad Tomovski, R. Hilfer, H. Srivastava,
\newblock{\it  Fractional and operational calculus with generalized fractional derivative operators and Mittag-Leffler type
functions, }
\newblock{Integral Transforms and Special Functions 21 (11) (2010)
797-814.}

\bibitem{JAN-DBD14}  J. A. Nanware, D. B. Dhaigude,
\newblock{\it Existence and uniqueness of solutions of differential equations of
fractional order with integral boundary conditions,}
\newblock{ J. Nonlinear Sci. Appl. 7, (2014), 246-254.}

\bibitem{RS-MS05}  R. K. Saxena, M. Saigo,
\newblock{\it Certain properties of fractional calculus
operators associated with generalized Mittag-Leffler function,}
\newblock{  Frac. Calc. Appl. Anal., Vol. 8 (2), (2005), 141-154. }

\bibitem{AS-DD18} Ahmad Y. A. Salamooni, D. D. Pawar,
\newblock{\it Existence and uniqueness of nonlocal boundary conditions
 for Hilfer-Hadamard-type fractional differential
equations.}
\newblock{ arXiv:1802.04262v1[math.AP] 12 Feb 2018. }



\end{thebibliography}
\end{document}